\documentclass[10pt]{amsart}

\font\de=cmssi10
\usepackage{amsmath}
\usepackage{amsthm}
\usepackage{amsfonts}

\newtheorem{teo}{\bf Theorem}

\newtheorem{prop}{\bf Proposition}

\begin{document}
\title[]{Hyperbolic sub-dynamics: compact invariant 3-manifolds}
\author{Jana Rodriguez Hertz}
\address{IMERL-Universidad de la Rep\'ublica}
\email{jana@fing.edu.uy}
\thanks{Partially supported by Fondo Clemente Estable 9021}
\begin{abstract}
In 1970, Hirsch asked what kind of compact invariant sets could be
part of a hyperbolic set. Here we obtain that, in case such an
invariant set is a 3D manifold, it is a connected sum of tori with
handles quotiented by involutions. Moreover, if the manifold is
orientable, the involutions are all trivial.\par In 1975, Ma{\~n}{\'e}
characterized hyperbolic dynamics restricted to manifolds and
called them quasi Anosov. We also classify here quasi-Anosov
dynamics in 3D-manifolds.
\end{abstract}
\maketitle \thispagestyle{empty} \vspace{-1.5em}
\section{Introduction}
In 1970, Hirsch asked what kind of compact invariant sets could
lie in a hyperbolic set. In particular, he asked whether the
restriction of a diffeomorphism to a hyperbolic set that is a
manifold induced an Anosov diffeomorphism, and found conditions
under which an affirmative answer is obtained \cite{hirsch}. In
1975, Ma\~n\'e found a characterization of hyperbolic dynamics
when restricted to compact invariant manifolds, which he called
{\de quasi Anosov diffeomorphisms} \cite{manhe}. Finally, in 1976,
Franks and Robinson gave an example of a quasi Anosov
diffeomorphism in a connected sum of two ${\mathbb T}^3$ that is
not Anosov \cite{frob}, giving a negative answer to the question
posed by Hirsch. This example consists essentially in considering
a linear Anosov system on a torus, say $T_1$, and its inverse on
another torus $T_2$. They produce appropriate perturbations on
each torus (DA diffeomorphisms) around their respective fixed
points. Then they cut suitable balls containing these fixed
points, and carefully glue together along their boundary so that
the stable and unstable foliations intersect quasi-transversally.
The aforementioned characterization by Ma{\~n}{\'e} yields a quasi-Anosov
diffeomorphism in the connected sum of $T_1$ and $T_2$, what
implies that $T_1\# T_2$ is a compact invariant subset of some
hyperbolic set.\par An example of a quasi Anosov diffeomorphism in
a non orientable 3-manifold may be found in \cite{zhumed}. The
example is similar to Franks and Robinson's, but a quotient by an
involution is done to the dynamics in the $T_i$'s before gluing
them together.
\par Here we show that {\em all} the examples of 3D compact invariant manifolds that lie in hyperbolic
sets are connected sums of these. More precisely,
\begin{teo}
Let $f:N\to N$ be a diffeomorphism, and let $M\subset N$ be a
hyperbolic set for $f$ such that $M$ is a 3D closed sub-manifold.
Then $M$ is the connected sum of a finite number tori with $r$
handles quotiented by involutions. That means, the Kneser
decomposition of $M$ is
$$M=\#_i({\mathbb T}^3/{\theta_i})\#r(S^1\times S^2)\qquad i>0,\quad r\geq 0$$
where all $\theta_i$ are involutions. \par Moreover, in case $M$
is orientable then all the involutions are trivial.
\end{teo}
We could also say that hyperbolic 3-manifolds are {\de generalized
connected sums} of $k$ tori quotiented by involutions, where the
generalized sum of $M_1$ and $M_2$ consists in taking away $r$
cells from $M_1$ and $M_2$ and gluing together along the boundary
of the cells by means of a reversing orientation
diffeomorphism.\par The dynamics on these sets is also classified:
\begin{teo}\label{teo.quasi.Anosov} Let $f$, $M$ and $N$ be as in theorem A.
And call $g=f|_M$ (so $g$ is a quasi Anosov diffeomorphism).
Then there exists $n>0$ such that
\begin{enumerate}
\item \label{item1} The non-wandering set $\Omega(g^n)$ consists
of $k$ connected codimension one expanding attractors $A_i$ and
$m$ connected codimension one shrinking repellors $R_j$. \item If
$\{W_l\}_{l=1}^{k+m}$ are the basins of $A_i$ and $R_j$, then each
restriction $g^n|{W_l}$ is topologically equivalent to quotient
functions $f_l^{n}/\theta_l$ acting on $({\mathbb T}^3\setminus
P_l)/\theta_l$ where
\begin{itemize}
\item $f_l$ are DA maps of ${\mathbb T}^3$,%
\item $P_l$ are finite sets of $f_l$-periodic points,%
\item and $\theta_l$ are involutions commuting with $f_l$
\end{itemize}   \item $M$ is the connected sum of $k+m$ (quotiented)
tori with handles. \item If $M$ is orientable then $\theta_l=id$
for all $l$.
\end{enumerate}
\end{teo}
Observe that all quasi Anosov diffeomorphisms can be obtained as
the restriction of a diffeomorphism to a hyperbolic set which is a
sub-manifold (see theorem \ref{teo.manhe}), so in particular we
obtain a classification of quasi-Anosov dynamics in dimension 3.
This completes a description started in \cite{rhuv}.\par Let us
remark that, from all examples above, only ${\mathbb T}^3$ can be
an invariant subset of any known Anosov system, due to a result by
A. Zeghib . In that case, the dynamics is Anosov. See
\cite{zeghib}\par Let us also mention that hyperbolic sets which
are 2D manifolds are always tori, and the dynamics on them is
Anosov. The question remains open
for hyperbolic submanifolds of higher dimensions .\par %
{\em Acknowledgemets} I thank the Department of Mathematics of the
U. of Toronto, where part of this work has been done, and M. Shub
for kind hospitality during my stay. I am also indebted to the
Departamento de Matem\'atica of Cuernavaca, specially A. Arroyo
and J. Seade. I thank F. Rodriguez Hertz, R. Ures and  A. Zeghib
for useful comments.
\section{Proof}
%
A diffeomorphism $g:M\to M$ is called {\de quasi Anosov} if
$\|Dg^n(x)v\|_n$ is unbounded for each non zero vector of $T_xM$.
There is a relation between quasi Anosov diffeomorphisms and {\de
hyperbolic sub-manifolds} (sub-manifolds that are hyperbolic
sets).
\begin{teo}\cite{manhe} \label{teo.manhe}The following statements are equivalent:
\begin{enumerate}
\item $g:M\to M$ is quasi Anosov \item There is a diffeomorphism
$f:N\to N$, and an embedding $i:M\hookrightarrow N$ verifying
$fi=ig$ such that $i(M)$ is a hyperbolic sub-manifold for $f$.
\item \label{item3} $g$ is Axiom A and all points $x\in M$ satisfy
$T_x W^s(x)\cap T_x W^u(x)=\{0\}$
\end{enumerate}
\end{teo}
Condition described in item (3) is called {\de
quasi-transversality}. Observe that, as a corollary, quasi Anosov
diffeomorphisms satisfy Axiom A and the no cycle condition. A
hyperbolic set $A$ is a {\de codimension one expanding attractor}
if $W^u(x)\subset A$ and $\dim W^u(x)=\dim M-1$ for all $x\in A$.
Codimension one shrinking repellors are defined analogously. A
codimension one expanding attractor is {\de orientable} if the
intersection index $W^s_{loc}(x)\cap W^u_{loc}(y)$ is constant at
all intersection points with  $x,y\in A$. Then we have the
following:
\begin{teo}\cite{ply}\label{teo.plykin}
If $A$ is a codimension one expanding attractor for $g:M^n\to M^n$
such that $W^s(A)$ is connected, then there is a DA map $f$ of
${\mathbb T}^n$, a finite set $P$ of $f$-periodic orbits, and an
involution $\theta$ of ${\mathbb T}^n$ conmuting with $f$, fixing
$P$, and making the following diagram commute:
$$\begin{array}{ccccc}
&             & f&               &\\
&{\mathbb T}^n\setminus P&\longrightarrow&{\mathbb T}^n\setminus P&\\
p_\theta&\downarrow&&\downarrow&p_\theta\\

&W^s(A)&\longrightarrow& W^s(A)&\\
&      & g &&
\end{array}\qquad
\begin{array}{l} \mbox{where}\quad p_\theta:{\mathbb T}^n\to {\mathbb T}^n/\theta\\ \mbox{is the canonical projection}
\end{array}
$$%

In particular, $W^s(A)$ is homeomorphic to ${\mathbb T}^n/\theta$
minus a finite number of points. If $A$ is orientable, then
$\theta=id$.
\end{teo}
Finally, let us state
\begin{teo}\cite{zhumed} If $M^3$ is orientable, then all codimension one expanding attractors are orientable.
\end{teo}
Since a codimension one expanding attractor of an Axiom A $g$ can
be written as a finite union of connected codimension one
expanding attractors of $g^n$ for some $n>0$, the whole thing is
reduced to proving $\Omega(g)$ consists of codimension one
attractors and repellors.
\begin{prop}
All attractors of 3D quasi Anosov diffeomorphisms are codimension one expanding attractors.
\end{prop}
Let us say that a basic set {\em is $(s,u)$} if its stable dimension is $s$ and its unstable dimension is $u$. Being
$\Lambda_i$ and $\Lambda_j$ basic sets, let us denote $\Lambda_i \to \Lambda_j$ if $W^u(\Lambda_i)\cap W^s(\Lambda_j)\ne \emptyset$.
Observe that this relation induces a connected graph in $\Omega(g)$ which has no cycles. Also, $\Lambda_i \to \Lambda_j$ implies that $s_i\geq s_j$ due to
item (\ref{item3}) in theorem \ref{teo.manhe} above.
\begin{prop}
If $A$ is an attractor and $\Lambda$ is a $(1,2)$ set such that $\Lambda\to A$ then $\Lambda=A$.
\end{prop}
Remark above immediately implies that $A$ is a codimension one
expanding attractor. Moreover, definition applies so $A\cup
W^u(\Lambda)$ is a codimension one expanding attractor. The rest
follows from Plykin's theorem above since $\Lambda \to A$ implies
$W^s(A\cup W^u(\Lambda))$ is connected.\par So, this proves that,
unless the diffeomorphism is Anosov, all maximal chains must have
a change of index. Also, that all maximal chains are of the form
$R\to A$ where $R$ is a codimension one shrinking repellor and $A$
is a codimension one expanding attractor. Each repellor/attractor
generates a 3-torus and, in case they are related, a connected sum
\`a la Franks-Robinson arises. Let us finally observe that these
models are not stable (indeed, it is easy to obtain perturbations
whose dynamic behavior is very different), but are
$\Omega$-stable.

\end{document}